\documentclass[a4paper,conference]{IEEEtran}
\IEEEoverridecommandlockouts

\usepackage{cite,color,url}
\usepackage{amsmath,amssymb,amsfonts,hyperref,multirow}
\usepackage{booktabs}
\usepackage{algorithmic,tikz}
\usepackage{subcaption,amsthm,longtable,lscape,amsmath}
\usepackage{graphicx,longtable}
\usepackage{textcomp}
\usepackage{xcolor}
\def\BibTeX{{\rm B\kern-.05em{\sc i\kern-.025em b}\kern-.08em
    T\kern-.1667em\lower.7ex\hbox{E}\kern-.125emX}}
\begin{document}

\title{On Solving the Shortest Paths with \\ Exclusive-Disjunction Arc Pairs Conflicts
}
\author{\IEEEauthorblockN{Roberto Montemanni}
\IEEEauthorblockA{\textit{Department of Sciences and Methods for Engineering} \\
\textit{University of Modena and Reggio Emilia}\\
Reggio Emilia, Italy \\
roberto.montemanni@unimore.it}
\and
\IEEEauthorblockN{Derek H. Smith}
\IEEEauthorblockA{\textit{Faculty of Computing, Engineering and Science} \\
\textit{University of South Wales}\\
Pontypridd, Wales, UK \\
derek.smith@southwales.ac.uk}}

\maketitle
\begin{abstract}
A variant of the well-known Shortest Path Problem is studied in this paper, where pairs of conflicting arcs are provided, and for each conflicting pair a penalty is paid once neither or both of the arcs are selected. This configures a set of soft-constraints. The problem, which can be used to model real applications, looks for a path from a given origin to a given destination that minimizes the cost of the arcs traversed plus the penalties incurred. 

In this paper, we consider a compact mixed integer linear
program representing the problem  and we solve it with the open-source solver CP-SAT, part of the Google OR-Tools computational suite.

An experimental campaign on the instances available from the literature indicates that the approach we propose achieves results comparable with those of  state-of-the-art solvers, notwithstanding it is a compact model, while the other approaches require the generation of dynamic constraints in order for the models to be competitive. Some best-known results have been improved in this study, and some instances have been closed for the first time.
\color{black}
\end{abstract}

\section{Introduction}
The shortest path problem is a classical optimization problem with several real-world applications, where the aim is to find the shortest (or quickest) way to go from one point to another on a given network \cite{meh08}. Several polynomial time algorithms exist for the problem \cite{dij59}. 

In this paper we consider the Shortest Paths with Exclusive-Disjunction Arc Pairs Conflicts (SP--EDAC), an extension of the basic problem, originally introduced in \cite{cer23}, where  conflicts exist between the arcs of some given pairs, and a penalty has to be paid if both or neither of the arcs of such pairs are part of the optimal solution.  The authors, on top of introducing the problem for the first time, propose two mixed integer linear programming models and a matheuristic approach based on one of the two models.

An application of the SP--EDAC problem is described in Cerulli et al. \cite{cer23}, in the context of software testing and validation. The nodes of the graph represent front-end pages of a web application. Every time a user interacts with a page, some server-side scripts are internally invoked in order to provide the requested services. The arcs of the graph model these actions, that as a side-effect take the user to other pages (nodes of the graph). Given a home page $s$ and a landing
page $t$, the aim is to find the shortest sequence of actions  to test the different server-side
scripts. In such a testing/validation context, conflicts are associated with actions activating the same script, and the aim is to test all the possible scripts no more than once during the journey from page $s$ to page $t$. 

Problems with elements in common with the SP--EDAC were previously treated in \cite{gab76}, where critical variants of the shortest path
problem arising in automatic generation of test paths for programs were introduced. A first definition of conflict pairs
was presented in \cite{kra73}, although nodes instead of arcs were part of the conflicts, and conflicts were hard constraints. A procedure to solve
the problem was discussed in \cite{sri82}.
A polyhedral study of the problem was presented in \cite{bla15}, while branch-and-bound and
dynamic programming approaches were recently proposed in \cite{fer21}.

In this paper, a  mixed integer linear programming model for the SP--EDAC is considered and solved by the open-source solver CP-SAT, which is part of the Google OR-Tools \cite{cpsat} optimization suite. Successful application of this solver on optimization problems with  characteristics similar to the problem under investigation, motivated our study \cite{md23}, \cite{cor}, \cite{rm25}. An experimental campaign on the benchmark instances previously used to validate algorithms in the literature is carried out.

The overall organization of the paper is as follows. The SP--EDAC problem is formally defined in Section \ref{desc}, while a  mixed integer linear program representing the problem is provided in Section \ref{model}. Section \ref{exp} is devoted to computational experiments: after the introduction of the instances commonly adopted in the literature, the results obtained by solving the model with the CP-SAT solver are compared with the published results on some of the instances. Conclusions about the work presented are finally drawn in Section \ref{conc}.

\section{Problem Description}\label{desc}

Given a directed graph $G = (V, A)$, with $V$ and $A$ being the sets
of vertices and arcs of $G$, respectively, the classical Shortest Path Problem can be defined as follows. A traversing cost $w_{ij} \in \mathbb{Z}^+_0$ is associated with each arc $(i,j) \in A$. Two vertices $s, t \in V$ are given as source and destination of the path.
The Shortest Path Problem seeks to find a path from $s$ to
$t$ that minimizes the sum of the costs of the arcs traversed. 

The SP--EDAC adds to the classical problem the concept of disjunctively exclusive conflicts on arc pairs:
if neither or both the arcs of a conflict pair  are selected, the conflict
is violated and a certain penalty has to be paid. Formally, a set $C$ of
conflicts is given: $C = \left \{ \{(i,j), (k,l)\}: (i,j), (k,l) \in A \text{ are in conflict}\right \}$ and a penalty $p_{ijkl}$ is associated to each $\{(i,j).(k,l)\} \in C$. The penalty is incurred once the arcs $(i,j)$ and $(k,l)$ are either both selected together or neither of them is. 
 
The SP--EDAC problem aims to minimize the overall cost of an $s-t$ path, given
by the sum of the costs of the selected arcs and the penalties of the
violated conflicts.

\begin{figure*}
{
\begin{center}
{
\begin{tikzpicture}[node distance={2.5cm}, main/.style = {draw, circle}]
			\node[main,minimum size=0.75cm] (f) {s};
			\node[main,minimum size=0.75cm] (a) [above right of=f] {a};
			\node[main,minimum size=0.75cm] (b) [below right of=f] {b};
			\node[main,minimum size=0.75cm] (c) [below right of=a] {c};
			\node[main,minimum size=0.75cm] (d) [above right of=c] {d};
			\node[main,minimum size=0.75cm] (e) [below right of=c] {e};
			\node[main,minimum size=0.75cm] (g) [below right of=d] {t};
			\draw [thick, line width=1.2,->,color=red] (f) to node[midway,above left] {3} (a) ;
			\draw [thick, line width=1.2,->,color=black] (f) to node[midway,below left] {1} (b) ;
			\draw [thick, line width=1.2,->,color=blue] (a) to node[midway,below left] {1} (c) ;
			\draw [thick, line width=1.2,->,color=green] (a) to node[midway,above] {2} (d) ;
			\draw [thick, line width=1.2,->] (b) to node[midway,left] {1} (a) ;
			\draw [thick, line width=1.2,->,color=black] (b) to node[midway,above left] {4} (c) ;
			\draw [thick, line width=1.2,->,color=blue] (b) to node[midway,below] {3} (e) ;
			\draw [thick, line width=1.2,->,color=black] (c) to node[midway,below right] {2} (d) ;
			\draw [thick, line width=1.2,->,color=black] (c) to node[midway,above right] {4} (e) ;
			\draw [thick, line width=1.2,->,color=red] (c) to node[midway,above] {2} (g) ;
			\draw [thick, line width=1.2,->,color=green] (d) to node[midway,above right] {1} (g) ;
			\draw [thick, line width=1.2,->,color=black] (e) to node[midway,below right] {3} (g) ;
\end{tikzpicture}
\hspace{1cm}
\begin{tikzpicture}[node distance={2.5cm}, main/.style = {draw, circle}]
			\node[main,minimum size=0.75cm] (f) {s};
			\node[main,minimum size=0.75cm] (a) [above right of=f] {a};
			\node[main,minimum size=0.75cm] (b) [below right of=f] {b};
			\node[main,minimum size=0.75cm] (c) [below right of=a] {c};
			\node[main,minimum size=0.75cm] (d) [above right of=c] {d};
			\node[main,minimum size=0.75cm] (e) [below right of=c] {e};
			\node[main,minimum size=0.75cm] (g) [below right of=d] {t};
			\draw [thick, line width=1.2,->,color=red] (f) to node[midway,above left] {3} (a) ;
			\draw [thick, line width=1.2,->,color=blue] (a) to node[midway,below left] {1} (c) ;
			\draw [thick, line width=1.2,->,color=black] (c) to node[midway,below right] {2} (d) ;
			\draw [thick, line width=1.2,->,color=green] (d) to node[midway,above right] {1} (g) ;
\end{tikzpicture}
}
	\caption{On the left an example of a small SP-EDAC instance is shown, where the  capacities are placed by the arcs. Arcs in black are not affected by conflicts while any two arcs sharing the same red, blue or green color are in conflict. On the right, an optimal solution that does not violate any conflict is shown. The zig-zag shape is motivated by the avoidance of not using a color or using the same color twice, in order to avoid any penalty. Notice that in larger examples it is common to incur some penalty in optimal solutions.}
	\label{figu}
\end{center}
}
\end{figure*}
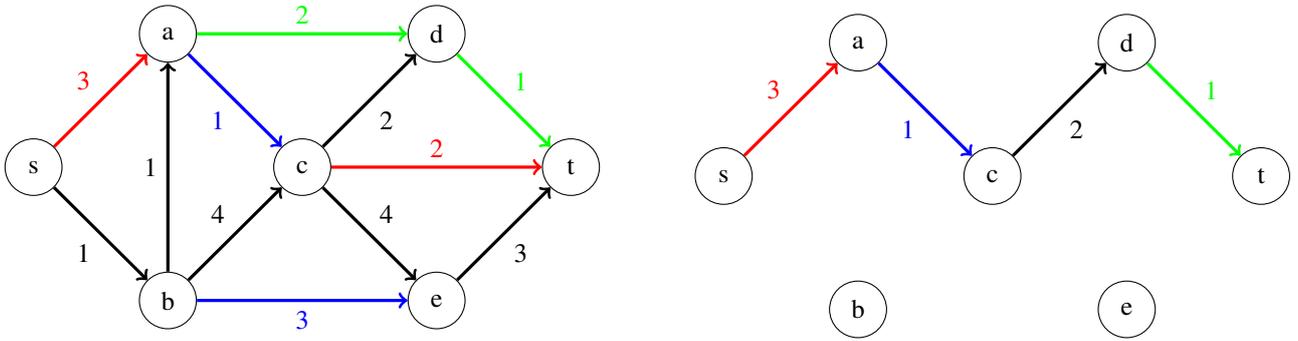

A small example of an SP-EDAC instance and an optimal solution for the instance are provided in Figure \ref{figu}.

\section{A Mixed Integer Linear Programming Model}\label{model}
In this section a model for the SP--EDAC is presented. It is inspired by the one originally proposed in \cite{cer23}, based on the linearization discussed in \cite{glo74}. The model is an extension of a classical model for the SP problem \cite{dij59}. Formally, a binary variable $x_{ij}$ is introduced for each $(i,j) \in A$. It takes value 1 if the arc $(i,j)$ is traversed by the path, 0 otherwise. Another set of variables is required to model the penalties of the SP--EDAC. A binary variable $y_{ijkl}$ is associated with a conflict $\{(i,j), (k,l)\} \in C$. The resulting  model is as follows.
\begin{align} 
\min & \sum_{(i,j) \in A} w_{ij}x_{ij} +  & \nonumber\\ 
 	& +  \!\!\!\!\! \sum_{\{(i,j),(k,l)\} \in C} \!\!\!\!\! p_{ijkl} (2y_{ijkl} - x_{ij} - x_{kl} +1)\!\!\!\!\!\!\!\!\!\!\!\!\!\!\!\!\!\!\!\!\!\!\!\!\!\!& \label{1}\\ 
s.t. &	\sum_{(j,i) \in A}  \!\!\!\! x_{ji}- \!\!\!\! \sum_{(i,j) \in A} \!\!\! x_{ij} \! = \! \begin{cases}
+1  \ \! \text{if $i=s$}\\
-1   \ \!\text{if $i=t$}\\
0   \ \ \ \!\text{if $i \notin \{s,t\}$}\\
\end{cases} \!\!\!\!\!\!\!\!\!\!\!\!\!\!\!\!\!\!\!\!\!\!\!\!\!\!\!\!\!\!\!\!\!\!\!\!\!\!\!\!\!\!\!\!\!\!\!\!\!\!\!\!\!\!\!\!\!\!\!\!\!\!\!\!\!\!\!\!\!\!\!\!\!\!\!\!\!\!\!\!& i \in V\label{2}\\
& \sum_{i,j \in S; (i,j) \in A} x_{ij} \le |S|-1 &  S \subseteq V, |S| \ge 2 \label{2a}\\
& y_{ijkl} \ge x_{ij} +x_{kl} -1 & \{(i,j), (k,l)\} \in C \label{3}\\
& y_{ijkl} \le x_{ij}   & \{(i,j), (k,l)\} \in C \label{4}\\
& y_{ijkl} \le x_{kl}  & \{(i,j), (k,l)\} \in C \label{4b}\\
& x_{ij} \in \{0,1\}& (i,j) \in A \label{5}\\
& y_{ijkl} \in \{0,1\} &  \{(i,j), (k,l)\} \in C \label{6}
\end{align}
The objective function (\ref{1}) minimizes the cost of the path and of the penalties incurred by such a path, in a form that produces  a stronger linear relaxation.
Constraints (\ref{2}) implement a classic flow model that guarantees feasible paths. 
Inequalities (\ref{2a}) are subtour elimination constraints \cite{dan54} that prevent the formation of loops created in order to avoid penalties.
Inequalities (\ref{3}) activate the penalty-variables once the two arcs of a conflicting pair are both in the solution.
Inequalites (\ref{4}) and (\ref{4b}) regulates the activation of the penalty-variables according to the objective function (\ref{1}).
Finally, the domains of the variables are specified in constraints (\ref{5}) and (\ref{6}).

\section{Computational Experiments} \label{exp}
Some implementation details behind the solution of the model discussed in Section \ref{model} with the solver selected are discussed in Section \ref{imp}.  The instances adopted for the experiments -- and commonly adopted in the previous literature of the SP--EDAC -- are described in Section \ref{ben}. An extensive comparison between the results achieved in this work and the results previously published, is presented in Section \ref{res}.

\subsection{Implementation Details}\label{imp}
In order to solve the model presented in Section \ref{model} to the Google OR-Tools CP-SAT solver \cite{cpsat} version 9.12, which does not support callbacks to alter the model during its solution, some adaptation are required to the model. The solver provides a construct called \emph{AddCircuit} that forces some of the arcs to conform like a loop touching a subset of the vertices of the graph. With the use of this construct it is possible to implement the model in a compact format. This construct is normally used to model generalized Traveling Salesman Problems \cite{mio} \cite{MC}, but it can be adapted to the purpose of modelling a path.

In order to implement the strategy, a new set of variables $x_{ii}$, which takes value 1 if vertex $i \in V \setminus \{s,t\}$ is not part of the optimal path, and 0 otherwise. Moreover, a constant $x_{ts}=1$ is added in order to close the optimal path from $s$ to $t$ in a loop.
Constraints (\ref{2}) and (\ref{2a}) are removed and at their place the following new constraints are inserted:
\begin{align} 
	& \text{AddCircuit}(x_{ij} | i, j \in V)&\label{7}\\
& x_{ii} \in \{0,1\}& i \in V \setminus \{s,t\}  \label{8}\\
&x_{ts}=1&\label{9}
\end{align}
Constraint (\ref{7}) imposes the presence of a circuit over the graph $G$. The definition of the variables in (\ref{8}) allows the circuit not to include some vertex $i \in V$ once the respective variables $x_{ii}$ takes value 1. The assignment of $x_{ts}$ to 1 finally impose the circuit to include the artificial arc $(t,s)$, imposing therefore that a path from $s$ to $t$ is in the circuit.

\subsection{Benchmark Instances}\label{ben}
Two sets of benchmark sets were proposed and made available in \cite{cer23} and represent the reference in the current literature.

The instances of the first set are purely random graphs, and an instances is generated according to the following parameters:
\begin{itemize}
\item $n$: number of nodes $|V|$, with \\ $n \in \{100, 200, 300, 400, 500\}$;
\item $d$: arc density, which defines the number of arcs as \\$|A|=d\cdot n(n-1)$,  with $d \in \{0.1, 0.2, 0.3\}$;
\item $r$: conflict density, which defines the number of conflicts as $|C| = \lfloor \frac{r}{2} \cdot m(m-1) \rfloor$, with $r \in \{10^{-3}, 2 \cdot 10^{-3}, 3 \cdot 10^{-3}\}$ when $n=100$, $r \in \{10^{-4}, 2 \cdot 10^{-4}, 3 \cdot 10^{-4}\}$ when $n=200$ and $r \in \{10^{-5}, 2 \cdot 10^{-5}, 3 \cdot 10^{-5}\}$ when $n \in \{300, 400, 500\}$;
\item $p$: penalty value range,  with \\$p \in \{[25, 125], [50, 150], [75, 175], [100, 200]\}$.
\end{itemize}
A total of 540 instances are present in the first benchmark set.

 The second benchmark set consists of graphs representing small-world 
networks. These networks are generated according \cite{wat98} and \cite{new99}). The instances are generated according to the following parameters:
\begin{itemize}
\item $n$: number of nodes $|V|$, with \\$n \in \{100, 200, 300, 400, 500\}$;
\item $k$: initial arc density. The number of neighbours for each node is initially $kn$, with \\$k \in \{0.15, 0.30, 0.45\}$;
\item $\beta$: rewiring probability: the probability that an arc is rewired, with $\beta = 0.5$;
\item $r$: conflict density, with $r \in \{10^{-3}, 2 \cdot 10^{-3}, 3 \cdot 10^{-3}\}$ when $n=100$, $r \in \{10^{-4}, 2 \cdot 10^{-4}, 3 \cdot 10^{-4}\}$ when $n=200$ and $r \in \{10^{-5}, 2 \cdot 10^{-5}, 3 \cdot 10^{-5}\}$ when $n \in \{300, 400, 500\}$;
\item $p$: penalty value range,  with $p \in \{[1, 20], [25, 200]\}$.
\end{itemize}
A total of 135 instances are present in the second benchmark set.

\subsection{Experimental Results} \label{res}
The experiments have been run on a computer equipped with an Intel Core i7 12700F CPU, while the results reported in \cite{cer23} were obtained on a computer based on an Intel Xeon E5 processor, with the latter being -- according to the conversion ratio inferred from \url{http://gene.disi.unitn.it/test/cpu_list.php} -- slightly slower. Moreover, we allow multithreading, while in \cite{cer23} the authors do not. This gives an  advantage to the approach we propose, although multithreading is at the basis of modern hardware architecture, and should probably be realistically exploited. A maximum computation time of 1800 seconds has been allowed for the approach we propose on each instance, against the 3600 seconds used by the methods proposed in \cite{cer23}. This should compensate for the difference of computing power while allowing all methods to reach steady results.

The results obtained by the approach we propose (CP-SAT) are compared with those of the three methods discussed in \cite{cer23}, namely a matheuristic approach and two Mixed Integer Linear Programming models (ILP and MILP), which are solved with IBM ILOG CPLEX 12.10.0 \cite{cplex}. In Tables \ref{t1} and \ref{t2}, and for each method considered, the average results on different subset of instances are reported, namely the value of the lower bounds (LB, when applicable), upper bounds (UB), time in seconds to retrieve the best solution (Sec Best) and total time in seconds (Sec Tot, when applicable) are provided. On top of this, for exact methods, the final percentage optimality gap (Opt gap \%) is also reported. The results are grouped by number of nodes $n$ and arc density ($d$ for random instances, $k$ for small-world instances). Cumulative averages over $d$ or $k$ are also reported. Whenever CP-SAT improves best-known results, the corresponding entries of the tables are highlighted in bold. 

The first observation about the results is that the matheuristic approach should be considered only if very fast results are required, otherwise the exact methods are able to improve its heuristic solutions, at the price of moderately longer computation times. When comparing exact methods, MILP seems to be dominant over ILP, as previously concluded in \cite{cer23}. The new solving model we propose, CP-SAT has performances comparable to those of MILP, although it seems to struggle more on the instances with small arc densities (especially on those with less nodes), and always keeping in mind that the experiments for CP-SAT have been run on a slightly faster machine. On the other hand, CP-SAT is able to improve the best-known results for some instance groups with a larger number of nodes and medium/high arc densities. Remarkably, all the small-world network instances with $k=0.30$  and all the random instances with $d=0.3$ and $n=400$ are closed here for the first time and the gap for the random instances with $d=0.2$ is substantially reduced. The average optimality gaps achieved by the exact methods are always negligible and below $0.5\%$, indicating that the available methods are already sufficiently good for the instances currently considered in the literature.

\begin{landscape}
\begin{table}[]
\caption{Experimental comparison on the random instances (average results).}\label{t1}
{
\scriptsize
\begin{tabular}{lrrrrrrrrrrrrrrrrr}
\toprule
                        & \multicolumn{2}{c}{Matheuristic \cite{cer23}}                  & \multicolumn{4}{c}{ILP \cite{cer23}}                                                                                           &                           & \multicolumn{5}{c}{MILP \cite{cer23}}                                                                                                                        & \multicolumn{5}{c}{CP-SAT}                                                                                                       \\
\cmidrule(lr){2-3}\cmidrule(lr){4-8}\cmidrule(lr){9-13}\cmidrule(lr){14-18}      
\multicolumn{1}{l}{Set} & \multicolumn{1}{r}{Cost} & \multicolumn{1}{r}{Sec} & \multicolumn{1}{r}{LB}      & \multicolumn{1}{r}{UB}      & \multicolumn{1}{r}{Sec} & \multicolumn{1}{r}{Sec}                               & \multicolumn{1}{r}{Opt} & \multicolumn{1}{r}{LB}       & \multicolumn{1}{r}{UB}       & \multicolumn{1}{r}{Sec} & \multicolumn{1}{r}{Sec}                                & \multicolumn{1}{r}{Opt} & \multicolumn{1}{r}{LB} & \multicolumn{1}{r}{UB} & \multicolumn{1}{r}{Sec} & \multicolumn{1}{r}{Sec}                             & \multicolumn{1}{r}{Opt} \\
\multicolumn{1}{r}{}    & \multicolumn{1}{r}{}     & \multicolumn{1}{r}{tot} & \multicolumn{1}{r}{}        & \multicolumn{1}{r}{}        & \multicolumn{1}{r}{best}  & \multicolumn{1}{r}{tot}   & \multicolumn{1}{r}{gap \%}      & \multicolumn{1}{r}{}         & \multicolumn{1}{r}{}         & \multicolumn{1}{r}{best}  & \multicolumn{1}{r}{tot}   & \multicolumn{1}{r}{gap \%}      & \multicolumn{1}{r}{}   & \multicolumn{1}{r}{}   & \multicolumn{1}{r}{best} & \multicolumn{1}{r}{tot} & \multicolumn{1}{r}{gap \%}      \\ 
\cmidrule(lr){1-1}\cmidrule(lr){2-3}\cmidrule(lr){4-8}\cmidrule(lr){9-13}\cmidrule(lr){14-18}      
$d=0.1, n=100$                & 72511.8                  & 22.6                    & \multicolumn{1}{r}{70649.1} & \multicolumn{1}{r}{70662.3} & \multicolumn{1}{r}{459.7} & \multicolumn{1}{r}{535.0} & 0.01871                   & {70658.9}  & {70658.9}  & \multicolumn{1}{r}{101.7} & \multicolumn{1}{r}{117.7} & 0.00000                   & 70522.1                & 70662.3                & 289.2                    & 765.0                 & 0.19832                   \\
$d=0.1, n=200$                & 125448.0                 & 8.0                     & 123950.3                    & 123950.3                    & 173.8                     & 202.7                     & 0.00000                   & \multicolumn{1}{r}{123950.3} & \multicolumn{1}{r}{123950.3} & \multicolumn{1}{r}{43.1}  & \multicolumn{1}{r}{46.4}  & 0.00000                   & 123950.3               & 123950.3               & 271.1                    & 399.0                   & 0.00000                   \\
$d=0.1, n=300$               & 50487.0                  & 7.4                     & 49850.1                     & 49850.1                     & 60.4                      & 61.8                      & 0.00000                   & 49850.1                      & 49850.1                      & 47.8                      & 52.3                      & 0.00000                   & 49850.1                & 49850.1                & 299.1                    & 416.5                   & 0.00000                   \\
$d=0.1, n=400$                 & 206286.7                 & 19.5                    & 205174.8                    & 205175.3                    & 410.0                     & 422.0                     & 0.00028                   & 205175.3                     & 205175.3                     & 231.6                     & 237.9                     & 0.00000                   & 205148.9               & 205180.9               & 891.9                    & 1112.9                  & 0.01560                   \\
$d=0.1, n=500$                 & 565514.5                 & 32.3                    & 563862.4                    & 563864.7                    & 707.4                     & 746.1                     & 0.00040                   & 563863.0                     & 563863.0                     & 395.3                     & 412.1                     & 0.00000                   & 563829.5               & 563872.5               & 1116.6                   & 1339.1                  & 0.00762                   \\
\midrule
 $d=0.1$                     & 204049.6                 & 18.0                    & 202697.3                    & 202700.6                    & 362.3                     & 393.5                     & 0.00158                   & 202699.5                     & 202699.5                     & 163.9                     & 173.3                     & 0.00000                   & 202660.2               & 202703.2               & 573.6                    & 806.5                 & 0.02122                   \\
\midrule
$d=0.2, n=100$                 & 362257.6                 & 16.9                    & 359479.1                    & 359480.7                    & 516.9                     & 570.0                     & 0.00044                   & 359480.7                     & 359480.7                     & 67.6                      & 77.1                      & 0.00000                   & 359465.9               & 359481.7               & 160.0                    & 288.7                   & 0.00440                   \\
$d=0.2, n=200$                 & 612722.6                 & 5.8                     & 611195.0                    & 611195.0                    & 125.6                     & 132.8                     & 0.00000                   & 611195.0                     & 611195.0                     & 48.9                      & 51.2                      & 0.00000                   & 611195.0               & 611195.0               & 189.3                    & 233.3                   & 0.00000                   \\
$d=0.2, n=300$                 & 290160.5                 & 9.4                     & 289097.6                    & 289097.6                    & 70.2                      & 71.2                      & 0.00000                   & 289097.6                     & 289097.6                     & 45.5                      & 46.9                      & 0.00000                   & 289097.6               & 289097.6               & 156.9                    & 188.0                   & 0.00000                   \\
$d=0.2, n=400$                 & 1007053.5                & 21.4                    & 1005089.6                   & 1005089.6                   & 463.7                     & 477.5                     & 0.00000                   & 1005089.6                    & 1005089.6                    & 178.6                     & 187.3                     & 0.00000                   & 1005084.2              & 1005089.6              & 445.7                    & 515.4                   & 0.00053                   \\
$d=0.2, n=500$                 & 2570130.5                & 38.6                    & 2567253.3                   & 2584767.1                   & 945.8                     & 1381.3                    & 0.67758                   & 2567254.7                    & 2577365.6                    & 483.1                     & 771.2                     & 0.39229                   & \textbf{2567259.9}     & \textbf{2567272.2}     & 599.5                    & 700.3                   & \textbf{0.00048}          \\
\midrule
$d=0.2$                     & 968465.0                 & 18.4                    & 966422.9                    & 969926.0                    & 424.4                     & 526.5                     & 0.36117                   & 966423.5                     & 968445.7                     & 164.8                     & 226.7                     & 0.20881                   & 966420.5               & \textbf{966427.2}      & 310.3                    & 385.1                   & \textbf{0.00069}          \\
\midrule
$d=0.3, n=100$                 & 874216.4                 & 20.4                    & 871319.1                    & 871350.5                    & 491.6                     & 634.5                     & 0.00360                   & 871340.3                     & 871340.3                     & 106.9                     & 115.0                     & 0.00000                   & 871284.1               & 871341.4               & 179.6                    & 360.7                   & 0.00658                   \\
$d=0.3, n=200$                 & 1461139.0                & 11.1                    & 1458897.4                   & 1458897.4                   & 199.5                     & 208.6                     & 0.00000                   & 1458897.4                    & 1458897.4                    & 57.6                      & 60.4                      & 0.00000                   & 1458897.4              & 1458897.4              & 186.6                    & 230.4                   & 0.00000                   \\
$d=0.3, n=300$                 & 716524.8                 & 13.2                    & 715210.8                    & 715210.8                    & 123.7                     & 127.3                     & 0.00000                   & 715210.8                     & 715210.8                     & 65.4                      & 67.3                      & 0.00000                   & 715210.8               & 715210.8               & 130.7                    & 153.3                   & 0.00000                   \\
$d=0.3, n=400$                 & 2387910.3                & 31.0                    & 2385794.9                   & 2385925.8                   & 675.8                     & 881.2                     & 0.00549                   & 2385795.9                    & 2385838.4                    & 438.9                     & 563.4                     & 0.00178                   & \textbf{2385804.0}     & \textbf{2385804.0}     & 225.9                    & 261.6                   & \textbf{0.00000}          \\
$d=0.3, n=500$                 & 5978248.7                & 75.4                    & 5975984.0                   & 5975986.0                   & 1270.3                    & 1321.8                    & 0.00003                   & 5975985.1                    & 5975985.1                    & 342.6                     & 350.3                     & 0.00000                   & 5975985.1              & 5975985.1              & 368.5                    & 426.6                   & 0.00000                   \\
\midrule
$d=0.3$                     & 2283607.8                & 30.2                    & 2281441.3                   & 2281474.1                   & 552.2                     & 634.7                     & 0.00144                   & 2281445.9                    & 2281454.4                    & 202.3                     & 231.3                     & 0.00037                   & 2281436.3              & \textbf{2281447.7}     & 218.2                    & 286.5                   & 0.00050                  \\
\bottomrule
\end{tabular}
}
\end{table}

\begin{table}[]
\caption{Experimental comparison on the small-world network instances (average results).}\label{t2}
{
\scriptsize
\begin{tabular}{lrrrrrrrrrrrrrrrrr}
\toprule
                        & \multicolumn{2}{c}{Matheuristic \cite{cer23}}                  & \multicolumn{4}{c}{ILP \cite{cer23}}                                                                                           &                           & \multicolumn{5}{c}{MILP \cite{cer23}}                                                                                                                        & \multicolumn{5}{c}{CP-SAT}                                                                                                       \\
\cmidrule(lr){2-3}\cmidrule(lr){4-8}\cmidrule(lr){9-13}\cmidrule(lr){14-18}      
\multicolumn{1}{l}{Set} & \multicolumn{1}{r}{Cost} & \multicolumn{1}{r}{Sec} & \multicolumn{1}{r}{LB}      & \multicolumn{1}{r}{UB}      & \multicolumn{1}{r}{Sec} & \multicolumn{1}{r}{Sec}                               & \multicolumn{1}{r}{Opt} & \multicolumn{1}{r}{LB}       & \multicolumn{1}{r}{UB}       & \multicolumn{1}{r}{Sec} & \multicolumn{1}{r}{Sec}                                & \multicolumn{1}{r}{Opt} & \multicolumn{1}{r}{LB} & \multicolumn{1}{r}{UB} & \multicolumn{1}{r}{Sec} & \multicolumn{1}{r}{Sec}                             & \multicolumn{1}{r}{Opt} \\
\multicolumn{1}{r}{}    & \multicolumn{1}{r}{}     & \multicolumn{1}{r}{tot} & \multicolumn{1}{r}{}        & \multicolumn{1}{r}{}        & \multicolumn{1}{r}{best}  & \multicolumn{1}{r}{tot}   & \multicolumn{1}{r}{gap \%}      & \multicolumn{1}{r}{}         & \multicolumn{1}{r}{}         & \multicolumn{1}{r}{best}  & \multicolumn{1}{r}{tot}   & \multicolumn{1}{r}{gap \%}      & \multicolumn{1}{r}{}   & \multicolumn{1}{r}{}   & \multicolumn{1}{r}{best} & \multicolumn{1}{r}{tot} & \multicolumn{1}{r}{gap \%}      \\ 
\cmidrule(lr){1-1}\cmidrule(lr){2-3}\cmidrule(lr){4-8}\cmidrule(lr){9-13}\cmidrule(lr){14-18}      
$k=0.15, n=100$                 & 82614.8                  & 25.1                    & 81041.8                & 81041.8                & 373.5                    & 477.2                   & 0.00000                   & 81041.8                & 81041.8                & 50.3                     & 68.2                    & 0.00000                   & 80949.2                & 81041.8                & 220.5                    & 602.4                   & 0.11421                   \\
$k=0.15, n=200$                & 166607.4                 & 7.1                     & 164923.3               & 164923.3               & 51.1                     & 53.0                    & 0.00000                   & 164923.3               & 164923.3               & 28.9                     & 29.9                    & 0.00000                   & 164923.3               & 164923.3               & 155.0                    & 234.8                   & 0.00000                   \\
$k=0.15, n=300$                 & 62683.3                  & 11.0                    & 61612.2                & 61612.2                & 84.1                     & 89.4                    & 0.00000                   & 61612.2                & 61612.2                & 69.7                     & 72.6                    & 0.00000                   & 61597.7                & 61612.2                & 504.0                    & 877.5                   & 0.02362                   \\
$k=0.15, n=400$                 & 268943.9                 & 25.6                    & 267468.8               & 267468.8               & 322.0                    & 333.2                   & 0.00000                   & 267468.8               & 267468.8               & 168.9                    & 200.3                   & 0.00000                   & 267432.9               & 267475.9               & 685.0                    & 888.9                   & 0.01608                   \\
$k=0.15, n=500$                 & 710501.0                 & 28.9                    & 708739.9               & 708739.9               & 387.9                    & 407.1                   & 0.00000                   & 708739.9               & 708739.9               & 265.9                    & 275.2                   & 0.00000                   & 708731.7               & 708749.4               & 786.1                    & 917.8                   & 0.00251                   \\
\midrule
$k=0.15$                     & 258270.1                 & 19.6                    & 256757.2               & 256757.2               & 243.7                    & 272.0                   & 0.00000                   & 256757.2               & 256757.2               & 116.7                    & 129.2                   & 0.00000                   & 256727.0               & 256760.5               & 470.1                    & 704.3                   & 0.01308                   \\
\midrule
$k=0.30, n=100$                 & 481112.0                 & 22.4                    & 478249.8               & 478257.8               & 580.7                    & 584.8                   & 0.00167                   & 478257.8               & 478257.8               & 49.3                     & 58.8                    & 0.00000                   & 478257.8               & 478257.8               & 116.6                    & 331.3                   & 0.00000                   \\
$k=0.30, n=200$                 & 791278.9                 & 7.7                     & 789284.8               & 789284.8               & 65.2                     & 69.8                    & 0.00000                   & 789284.8               & 789284.8               & 48.7                     & 49.9                    & 0.00000                   & 789284.8               & 789284.8               & 73.6                     & 108.5                   & 0.00000                   \\
$k=0.30, n=300$                 & 376861.7                 & 9.9                     & 375456.3               & 375456.3               & 81.1                     & 83.8                    & 0.00000                   & 375456.3               & 375456.3               & 55.0                     & 57.9                    & 0.00000                   & 375456.3               & 375456.3               & 179.9                    & 223.2                   & 0.00000                   \\
$k=0.30, n=400$                 & 1294249.8                & 37.0                    & 1292325.1              & 1292325.1              & 260.5                    & 266.6                   & 0.00000                   & 1292325.1              & 1292325.1              & 163.8                    & 166.8                   & 0.00000                   & 1292325.1              & 1292325.1              & 325.5                    & 370.0                   & 0.00000                   \\
$k=0.30, n=500$                 & 3293128.1                & 59.6                    & 3290971.6              & 3291039.4              & 1330.3                   & 1434.2                  & 0.00206                   & 3290979.9              & 3290982.9              & 605.1                    & 716.4                   & 0.00009                   & \textbf{3290981.3}     & \textbf{3290981.3}     & 472.7                    & 549.7                   & \textbf{0.00000}          \\
\midrule
$k=0.30$                     & 1247326.1                & 27.3                    & 1245257.5              & 1245272.7              & 463.6                    & 487.8                   & 0.00122                   & 1245260.8              & 1245261.4              & 184.4                    & 210.0                   & 0.00005                   & \textbf{1245261.1}     & \textbf{1245261.1}     & 233.6                    & 316.5                   & \textbf{0.00000}          \\
\midrule
$k=0.45, n=100$                 & 1088992.7                & 18.1                    & 1085908.1              & 1085908.1              & 452.3                    & 529.0                   & 0.00000                   & 1085908.1              & 1085908.1              & 48.8                     & 51.2                    & 0.00000                   & 1085867.2              & 1085922.1              & 178.3                    & 317.6                   & 0.00505                   \\
$k=0.45, n=200$                 & 1880620.4                & 11.5                    & 1878991.4              & 1878991.4              & 172.5                    & 198.8                   & 0.00000                   & 1878991.4              & 1878991.4              & 34.7                     & 37.8                    & 0.00000                   & 1878991.4              & 1878991.4              & 134.1                    & 151.6                   & 0.00000                   \\
$k=0.45, n=300$                 & 909540.1                 & 15.6                    & 907954.4               & 907954.4               & 110.1                    & 111.3                   & 0.00000                   & 907954.4               & 907954.4               & 39.0                     & 40.3                    & 0.00000                   & 907954.4               & 907954.4               & 113.0                    & 132.5                   & 0.00000                   \\
$k=0.45, n=400$                 & 3069898.0                & 52.1                    & 3067816.6              & 3067816.6              & 299.9                    & 307.1                   & 0.00000                   & 3067816.6              & 3067816.6              & 136.2                    & 138.0                   & 0.00000                   & 3067816.6              & 3067816.6              & 301.1                    & 336.7                   & 0.00000                   \\
$k=0.45, n=500$                & 7544776.6                & 88.9                    & 7541803.7              & 7541803.7              & 1155.0                   & 1171.7                  & 0.00000                   & 7541803.7              & 7541803.7              & 237.2                    & 246.4                   & 0.00000                   & 7541803.7              & 7541803.7              & 385.0                    & 431.8                   & 0.00000                   \\
\midrule
$k=0.45$                     & 2898765.6                & 37.2                    & 2896494.8              & 2896494.8              & 438.0                    & 463.6                   & 0.00000                   & 2896494.8              & 2896494.8              & 99.2                     & 102.7                   & 0.00000                   & 2896486.7              & 2896497.6              & 222.3                    & 274.0                   & 0.00038                  \\
\bottomrule
\end{tabular}
}
\end{table}
\end{landscape}

In conclusion, modelling the SP--EDAC problem in terms of CP-SAT and attacking it with such a solver is a viable and practical way, especially when considering that such an approach rules out all the implementation burden related to the use of subtour elimination constraints otherwise required \cite{cer23}.
\color{black}

\section{Conclusions} \label{conc}
A  formulation based on mixed integer linear formulation for the Shortest Paths with Exclusive-Disjunction Arc Pairs Conflicts was considered and solved via the open-source solver CP-SAT, part of the Google OR-Tools computational suite. 

The experimental results indicate that the approach we propose has results comparable with those of the state-of-the-art solvers available in the literature, notwithstanding it does not required the complex programming techniques required by the other methods. Remarkably some of the best-known results have been improved by the approach we proposed.

\section*{Acknowledgments}
The work was partially supported by the Google Cloud Research Credits Program.
\color{black}
\bibliographystyle{IEEEtran}
\bibliography{mybibfile}

\end{document}